\documentclass[11pt]{article}
\usepackage{amssymb,amsmath,amsthm,pslatex,times,subfigure,url}

\title{Uncoverings on graphs and network reliability}
\author{Robert F.~Bailey\footnote{
Department of Mathematics and Statistics, University of Regina, 3737 Wascana Parkway, Regina, Saskatchewan S4S 0A2, Canada.  E-mail: \texttt{robert.bailey@uregina.ca}}\hspace{2mm} and
Brett Stevens\footnote{School of Mathematics and Statistics, Carleton University, 1125 Colonel By Drive, Ottawa, Ontario K1S 5B6, Canada.  E-mail: \texttt{brett@math.carleton.ca}}}

\newtheorem{thm}{Theorem}
\newtheorem{lemma}[thm]{Lemma}
\newtheorem{prop}[thm]{Proposition}

\newtheorem*{conj}{Conjecture}

\theoremstyle{definition}
\newtheorem{defn}[thm]{Definition}
\newtheorem{example}[thm]{Example}
\newtheorem{construction}[thm]{Construction}

\newtheorem*{proofcorr}{Proof of correctness}
\newtheorem*{proofmin}{Proof of minimality}

\renewcommand{\emptyset}{\varnothing}

\newcommand{\blob}{\circle*{0.2}}

\begin{document}

\maketitle

\begin{abstract}
We propose a network protocol similar to the $k$-tree protocol of Itai and Rodeh [{\em Inform.\ and Comput.}\ {\bf 79} (1988), 43--59].  To do this, we define a {\em $t$-uncovering-by-bases} for a connected graph $G$ to be a collection $\mathcal{U}$ of spanning trees for $G$ such that any $t$-subset of edges of $G$ is disjoint from at least one tree in $\mathcal{U}$, where $t$ is some integer strictly less than the edge connectivity of $G$.  We construct examples of these for some infinite families of graphs.  Many of these infinite families utilise factorisations or decompositions of graphs. In every case the size of the uncovering-by-bases is no larger than the number of edges in the graph and we conjecture that this may be true in general.
\end{abstract}



\section{Introduction}

In \cite{Itai88}, Itai and Rodeh proposed a communication protocol, which we shall call the {\em $t$-tree protocol}, which allows all nodes of a network to communicate through a distinguished root node, even when some set of $t-1$ or fewer edges are removed from the network.  
(In their paper, they used $k$ rather than $t$; we have changed the notation to be consistent with the conventions of design theory, because of the connections to that subject in this paper.)  
The protocol requires the graph $G$ modelling the network to have two properties.  First, the graph, $G$, must remain connected when any $t-1$ edges are removed.
Second, for any vertex, $r$, it requires a collection of $t$ spanning trees for $G$, $\{T_1,\ldots,T_t\}$, with the following property (the {\em $t$-tree condition for edges}): for any vertex $v$ and any $i$, $j$, where $1 \leq i < j \leq t$, the paths in $T_i$ and $T_j$ from $v$ to $r$ are internally disjoint.  A set of $t$ disjoint spanning trees clearly satisfies this condition and so can be used in the protocol for robust communication \cite{maxstp}.

Of course, for a given $t$, an arbitrary graph $G$ may not meet the $t$-tree condition.  An alternative solution for robust communications is to use a collection of $s$ (which may be greater than $t$) spanning trees (which need not necessarily be edge-disjoint) spanning trees, such that if any $t-1$ edges of $G$ are removed, at least one of the spanning trees remains intact.  Ideally, we would want this collection to be as small as possible (or of bounded size), and for the value of $t$ to be as large as possible.  The purpose of this article is to study such collections of spanning trees.  We begin with some basic definitions.

An {\em edge cut} in $G$ is a partition $(V_1,V_2)$ of the vertex set of $G$ into two non-empty subsets.  In other words, an edge cut is a set where the removal of the edges between $V_1$ and $V_2$ disconnects $G$; if the number of such edges is $t$, we call it a {\em $t$-edge cut}.  The {\em edge connectivity} of $G$ is the least value of $t$ for which there exists a $t$-edge cut in $G$; we denote this by $\lambda(G)$.  We say that $G$ is {\em $t$-edge connected} if $\lambda(G) \geq t$.  We note that, by abuse of terminology, we will sometimes refer to a $t$-edge cut by the set of edges whose removal disconnects the graph, rather than the partition of $V$.

\begin{defn} \label{defn:uncovering}
Let $n$, $k$ and $t$ be positive integers satisfying $n>k$ and $t\leq n-k$, and let $X$ be a set of size $n$.  An \emph{$(n,k,t)$-uncovering} is a collection $\mathcal{U}$ of $k$-subsets of $X$ such that any $t$-subset of $X$ is disjoint from at least one $k$-subset in $\mathcal{U}$.
\end{defn}

An $(n,k,t)$-uncovering is equivalent to an \emph{$(n,n-k,t)$ covering design}, which is a set of $(n-k)$-subsets, called blocks, such that any $t$-subset is contained in at least one block.  So results on coverings, such as those in the survey by Mills and Mullin \cite{MillsMullin92}, also give us results on uncoverings.  In both cases, the problem is to find an (un)covering of least possible size.  The most general bound is the {\em Sch\"{o}nheim bound}, which gives a lower bound of
\[ \left\lceil \frac{n}{n-k} \left\lceil \frac{n-1}{n-k-1} \left\lceil \cdots \left\lceil \frac{n-t+1}{n-k-t+1} \right\rceil \cdots \right\rceil \right\rceil \right\rceil \]
on the size of such an $(n,k,t)$-uncovering.  Uncoverings were introduced and studied by the first author in \cite{thesis,btubb,ecpg}; they were later introduced independently by Kroll and Vincenti \cite{Kroll08} by the name {\em antiblocking systems}.

In this paper, we are concerned with the case where $X$ is the edge-set of a graph $G$, and where each member of the uncovering is a spanning tree.  That is, we have the following.

\begin{defn} \label{defn:UBB}
Let $G=(V,E)$ be a connected graph and $t$ a positive integer.  A {\em $t$-uncovering-by-bases} for $G$ is a collection $\mathcal{U}$ of spanning trees for $G$ such that any $t$-subset of $E$ is disjoint from at least one spanning tree in $\mathcal{U}$.
\end{defn}

The name ``uncovering-by-bases'' (or UBB for short) comes from the fact that the spanning trees of a connected graph $G$ are precisely the bases of the cycle matroid $M(G)$ (see Oxley \cite{Oxley92} for further details).  It is possible to define UBBs for arbitrary matroids, as is done in \cite{btubb}.
Note that for a $t$-UBB to exist, we require that the graph obtained by deleting an arbitrary $t$-subset of edges from $G$ must have a spanning tree, which happens if and only if it is connected.  In other words, we require that the edge connectivity, $\lambda(G)$, must be strictly greater than $t$, so $t \leq \lambda(G)-1$.  In this paper, we only consider the case where this maximum is achieved, i.e.\ when $t=\lambda(G)-1$.

If $G$ happens to have $t+1=\lambda(G)$ edge-disjoint spanning trees, then we can use these as a $t$-UBB.  We call such graphs {\em maximum spanning tree-packable}, or {\em max-STP} graphs; these are described in \cite{maxstp}.  In a max-STP graph $G$, the $t$-UBB formed of the collection of $\lambda(G)$ edge-disjoint spanning trees is therefore optimal in two ways: first, because the number of edges $t$ which can be uncovered is as large as possible; second, because the spanning trees are all edge-disjoint, the size of the $t$-UBB is as small as possible.

\section{A collection of examples}
In this section, we present some constructions of UBBs for certain families or classes of graphs.  In each case, our constructions have the highest ``uncovering'' ability: we are able to take $t=\lambda(G)-1$.

\subsection{Complete bipartite graphs} \label{subsect:Kmn}
Consider the complete bipartite graph $K_{m,n}$.  Suppose that $2 \leq m \leq n$, and that the vertex set is $X \dot{\cup} Y$, where $|X|=m$ and $|Y|=n$.  Now, the edge connectivity of $K_{m,n}$ is $\min\{m,n\}=m$.  We show how to construct a $t$-UBB for $K_{m,n}$ where $t=m-1$ (i.e.~the largest possible t).

\begin{figure}[hbtp]
\centering
\setlength{\unitlength}{7mm}
\begin{picture}(5,6)
\put(1,6){\circle*{0.2}}
\put(1,5){\circle*{0.2}}
\put(1,4){\circle*{0.2}}
\put(1,2.3){\vdots}
\put(1,1){\circle*{0.2}}

\put(4,6){\circle*{0.2}}
\put(4,5){\circle*{0.2}}
\put(4,4){\circle*{0.2}}
\put(4,2.3){\vdots}
\put(4,1){\circle*{0.2}}
\put(4,0){\circle*{0.2}}

\put(1,6){\line(1,0){3}}
\put(1,6){\line(3,-1){3}}
\put(1,6){\line(3,-2){3}}
\put(1,6){\line(3,-5){3}}
\put(1,6){\line(1,-2){3}}

\put(4,6){\line(-3,-1){3}}
\put(4,6){\line(-3,-2){3}}
\put(4,6){\line(-3,-5){3}}

\put(0.5,5.85){$u$}
\put(4.2,5.85){$v$}
\end{picture}
\caption{\label{fig:Kmn} An example of a spanning tree $T_{uv}$ for $K_{m,n}$.}
\end{figure}
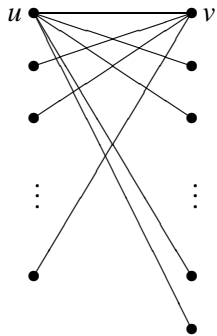

Fix an $m$-subset $S\subseteq Y$.  Let $A$ be an arbitrary $t$-subset of edges.  Now, since $|A|<m$, there exists a vertex $u \in X$ incident with no edge of $A$.  Similarly, since $|A|<|S|=m$, there exists $v \in S$ incident with no edge of $A$.  Now we construct a spanning tree $T_{uv}$ which contains the edge $uv$, and all other edges incident with each of $u$ and $v$. By construction, the set of spanning trees $\mathcal{U} = \{T_{uv} \,\mid\, u\in X, v\in S \}$ is a $t$-UBB for $K_{m,n}$, of size $m^2$.  A typical spanning tree $T_{uv}$ is shown in Figure \ref{fig:Kmn}.

We remark that this construction gives a minimal UBB: this can be seen by considering a matching in $K_{m,n}$ of size $m-1$.  However, we do not claim that it is minimum.  We also notice that the number of spanning trees in $\mathcal{U}$ is bounded above by the number of edges of $K_{m,n}$; the bound is sharp, as it is achieved when $m=n$.

\subsection{Graphs with Hamiltonian decompositions} \label{subsect:hamdec}
We now consider graphs with {\em Hamiltonian decompositions}, i.e.\ graphs which admit a partition of the edge set into Hamilton cycles.  In particular, we note that if $n$ is odd, the complete graph $K_n$ admits a Hamiltonian decomposition; this was known to Walecki in the 1890s (see Bryant \cite{Bryant07} for details).  A survey of more general results about Hamiltonian decompositions can be found in Section 4 of Gould \cite{Gould03}.  

To construct uncoverings-by-bases for these graphs, we need to know their edge connectivity.

\begin{prop} \label{prop:EdgeConnHamDec}
Let $G$ be a graph with a Hamiltonian decomposition into $c$ cycles.  Then the edge connectivity of $G$ is $\lambda(G)=2c$.
\end{prop}

\begin{proof}  Because $G$ has a Hamiltonian decomposition into $c$ cycles, $G$ must be $2c$-regular, so $\lambda(G) \leq 2c$.  Also, any edge-cut of $G$ must contain at least two edges from each of the $c$ Hamilton cycles, so $\lambda(G) \geq 2c$.  Hence $\lambda(G)=2c$.  \end{proof}

The construction works as follows.  

\begin{construction} \label{constr:HamDec}
Let $\mathcal{D} = \{C_1,\ldots,C_c\}$ be a Hamiltonian decomposition of $G$.  For each $C_i\in \mathcal{D}$ and for each $e \in C_i$, form a path $C_i\setminus e$, which is a spanning tree.  We claim that the set of all such paths is a $t$-UBB for $G$, where $t=2c-1$.
\end{construction}

\begin{proofcorr}
Let $A$ be an arbitrary $t$-subset of $E(G)$.  Now, if there exists a cycle $C_i$ such that $C_i \cap A = \emptyset$, then any path in $C_i$ is disjoint from $A$. So we suppose not, i.e.\ we suppose that each of the $c$ cycles contains some of the $t=2c-1$ edges in $A$.
By counting, there must exist a cycle $C_j$ containing exactly one edge $e\in A$, so the path $C_j\setminus e$ is disjoint from $A$. \qed
\end{proofcorr}

By construction, this UBB is minimal: if the $2c-1$ ``bad'' edges are arranged so that there is one bad edge $e$ in $C_1$ and two bad edges in each of $C_2,\ldots,C_c$, then the only spanning tree avoiding these edges is $C_1\setminus e$.  As in the previous subsection, however, we don't claim that this UBB is minimum.

\begin{figure}[hbtp]
\setlength{\unitlength}{7mm}
\centering
\begin{picture}(6,6)
\put(1.8,0.3){\line(1,0){2.25}}
\put(0.3,1.8){\line(0,1){2.25}}
\put(5.55,1.8){\line(0,1){2.25}}

\put(4.05,0.3){\line(1,1){1.5}}
\put(1.8,0.3){\line(-1,1){1.5}}

\put(4.05,0.3){\line(2,5){1.5}}

\put(1.8,0.3){\line(5,2){3.75}}
\put(1.8,0.3){\line(-2,5){1.5}}
\put(0.3,1.8){\line(5,-2){3.75}}

\put(0.3,4.05){\line(1,0){5.25}}

\put(1.8,0.3){\circle*{0.2}}
\put(4.05,0.3){\circle*{0.2}}
\put(5.55,1.8){\circle*{0.2}}
\put(5.55,4.05){\circle*{0.2}}
\put(0.3,4.05){\circle*{0.2}}
\put(0.3,1.8){\circle*{0.2}}
\put(2.925,5.55){\circle*{0.2}}

\qbezier(0.3,4.05)(1.6125,4.8)(2.925,5.55)
\qbezier(5.55,4.05)(4.2375,4.8)(2.925,5.55)

\qbezier(0.3,1.8)(1.6125,3.675)(2.925,5.55)
\qbezier(5.55,1.8)(4.2375,3.675)(2.925,5.55)

\end{picture}
\caption{A circulant graph on 7 vertices. \label{fig:circulant}}
\end{figure}
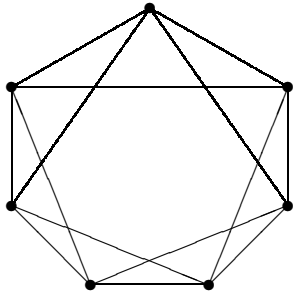

\begin{example}
Consider the graph shown in Figure \ref{fig:circulant} (an example of a {\em circulant graph}).  This is the union of two Hamilton cycles of length 7, so is 4-edge connected and we have $t=3$.  In each cycle we obtain 7 paths, and a total of 14 paths in our UBB.
\end{example}

We note that the number of spanning trees in the $t$-UBB is precisely the number of edges of $G$; so as in the previous subsection $|E(G)|$ is an upper bound on the minimum size of a $t$-UBB.  We also remark that this construction is very similar to the construction of a UBB for the permutation group $S_n$ in its action on 2-subsets: see \cite{ecpg} for details.

\subsection{Using 1-factorisations of graphs} \label{subsect:HKL}
In the previous subsection, we were able to construct UBBs for complete graphs with an odd number of vertices by virtue of the fact that they have Hamiltonian decompositions.  In order to consider a class of graphs which includes complete graphs with an even number of vertices, we will consider 1-factorisations with some particularly useful properties.

A {\em 1-factor} (also known as a {\em perfect matching}) in a graph $G$ is a spanning, 1-regular subgraph of $G$ (i.e.~a collection of vertex-disjoint edges of $G$ incident with every vertex).  A {\em 1-factorisation} of $G$ is a partition of the edge set of $G$ into edge-disjoint 1-factors.  A {\em perfect 1-factorisation} of $G$ is a 1-factorisation $\mathcal{F} = \{F_1,\ldots,F_k\}$ of $G$ where 
the union of any two 1-factors is a Hamilton cycle; 
if $G$ possesses a perfect 1-factorisation it is called {\em strongly Hamiltonian}.  (See Andersen \cite{Andersen07} for more detail on these topics.)  This notion of strong Hamiltonicity was introduced by Kotzig and Labelle in 1978 \cite{Kotzig78}, where they consider in detail cubic (i.e.~3-regular) strongly Hamiltonian graphs.

In general, it is difficult to determine if a graph has a perfect 1-factorisation.  It is conjectured, but not known, that the complete graph $K_{2m}$ on an even number of vertices $2m$ always has such a 1-factorisation (see Andersen \cite{Andersen07} or Wanless \cite{Wanless05}).  However, this property is often much stronger than we require, and so we introduce the following idea.  To every 1-factorisation $\mathcal{F}$ of a graph $G$, we can build an auxiliary digraph, $\mathcal{H}(G,\mathcal{F})$, with the 1-factors as vertices and an arc between two 1-factors $F$ and $F'$ if their union is a Hamilton cycle.
We note that whenever $(F,F')$ is an arc, then so is $(F',F)$.  We have chosen to express the adjacences by directed edges for two reasons: first, so that 2-cycles are well defined; second, so that we have a good notion of the successor of a vertex along an arc.

We are especially interested in graphs $G$ with $n$ vertices and which have a 1-factorisation $\mathcal{F}$ where the auxiliary digraph $\mathcal{H}(G,\mathcal{F})$ has a directed 2-factor (i.e.\ a spanning subgraph formed of directed cycles).  If there is a directed 2-factor which is a Hamilton cycle, then the 1-factorisation $\mathcal{F}$ is said to be {\em sequentially perfect}, as studied by Dinitz, Dukes and Stinson \cite{MR2134164}.  It is easy to see that for a perfect 1-factorisation $\mathcal{F}$, the digraph $\mathcal{H}(G,\mathcal{F})$ is a complete digraph, and thus $\mathcal{F}$ is sequentially perfect.

Also, we define an {\em HKL decomposition}\footnote{
The initials HKL were chosen in honour of W.~R.~Hamilton, A.~Kotzig and J.~Labelle.%
} of a graph $G$ to be a partition of the edges of $G$ into Hamilton cycles (on $n$ vertices) and a cubic strongly Hamiltonian graph (also on $n$ vertices).  In other words $\mathcal{H}(G,\mathcal{F})$ has a 2-factor with all components of size 2 except for one of size 3.  Such decompositions arise in the context of random regular graphs, which we will discuss in Section~\ref{random_regular}.  The following is an example of a graph of odd degree with an HKL decomposition.

\begin{example} \label{example:HKLexample}
Consider the graph in Figure \ref{figure:HKL1}(a), which has an HKL decomposition into a single Hamilton cycle and a cubic strongly Hamiltonian graph, as shown in Figure \ref{figure:HKL1}(b).

\begin{figure}[hbtp]
\setlength{\unitlength}{7mm}
\centering

\subfigure[A graph with an HKL decomposition]{%
\begin{picture}(10,5)

\put(4.000,0.000){\blob}
\put(2.586,1.414){\blob}
\put(2.586,3.414){\blob}
\put(4.000,4.828){\blob}
\put(6.000,4.828){\blob}
\put(7.414,3.414){\blob}
\put(7.414,1.414){\blob}
\put(6.000,0.000){\blob}

\put(2.586,3.414){\line(1,1){1.414}}
\put(4.000,4.828){\line(1,0){2.000}}
\put(6.000,4.828){\line(1,-1){1.414}}
\put(7.414,1.414){\line(0,1){2.000}}
\put(6.000,0.000){\line(1,1){1.414}}
\put(4.000,0.000){\line(1,0){2.000}}
\put(2.586,1.414){\line(1,-1){1.414}}
\put(2.586,1.414){\line(0,1){2.000}}

\qbezier(2.586,1.414)(3.293,3.121)(4.000,4.828)
\qbezier(7.414,1.414)(6.707,3.121)(6.000,4.828)
\put(4.000,0.000){\line(1,1){3.414}}
\put(6.000,0.000){\line(-1,1){3.414}}

\put(4.000,0.000){\line(0,1){4.828}}
\put(6.000,0.000){\line(0,1){4.828}}
\put(2.586,1.414){\line(1,0){4.828}}
\put(4.000,4.828){\line(1,-1){3.414}}

\qbezier(4.000,0.000)(3.293,1.707)(2.586,3.414)
\qbezier(6.000,0.000)(6.707,1.707)(7.414,3.414)
\qbezier(2.586,3.414)(4.293,4.121)(6.000,4.828)
\qbezier(2.586,1.414)(5.000,2.414)(7.414,3.414)

\end{picture}  }

\subfigure[The decomposition into a cubic strongly Hamiltonian graph (left) and a Hamilton cycle (right).]{%

\begin{picture}(6,5)

\put(2.000,0.000){\blob}
\put(0.586,1.414){\blob}
\put(0.586,3.414){\blob}
\put(2.000,4.828){\blob}
\put(4.000,4.828){\blob}
\put(5.414,3.414){\blob}
\put(5.414,1.414){\blob}
\put(4.000,0.000){\blob}

\put(0.586,3.414){\line(1,1){1.414}}
\put(2.000,4.828){\line(1,0){2.000}}
\put(4.000,4.828){\line(1,-1){1.414}}
\put(5.414,1.414){\line(0,1){2.000}}
\put(4.000,0.000){\line(1,1){1.414}}
\put(2.000,0.000){\line(1,0){2.000}}
\put(0.586,1.414){\line(1,-1){1.414}}
\put(0.586,1.414){\line(0,1){2.000}}

\qbezier(0.586,1.414)(1.293,3.121)(2.000,4.828)
\qbezier(5.414,1.414)(4.707,3.121)(4.000,4.828)
\put(2.000,0.000){\line(1,1){3.414}}
\put(4.000,0.000){\line(-1,1){3.414}}

\end{picture}
\qquad
\begin{picture}(6,5)

\put(2.000,0.000){\blob}
\put(0.586,1.414){\blob}
\put(0.586,3.414){\blob}
\put(2.000,4.828){\blob}
\put(4.000,4.828){\blob}
\put(5.414,3.414){\blob}
\put(5.414,1.414){\blob}
\put(4.000,0.000){\blob}
\put(2.000,0.000){\line(0,1){4.828}}
\put(4.000,0.000){\line(0,1){4.828}}
\put(0.586,1.414){\line(1,0){4.828}}
\put(2.000,4.828){\line(1,-1){3.414}}

\qbezier(2.000,0.000)(1.293,1.707)(0.586,3.414)
\qbezier(4.000,0.000)(4.707,1.707)(5.414,3.414)
\qbezier(0.586,3.414)(2.293,4.121)(4.000,4.828)
\qbezier(0.586,1.414)(3.000,2.414)(5.414,3.414)

\end{picture} }
\caption{An example of an HKL decomposition. \label{figure:HKL1}}
\end{figure}
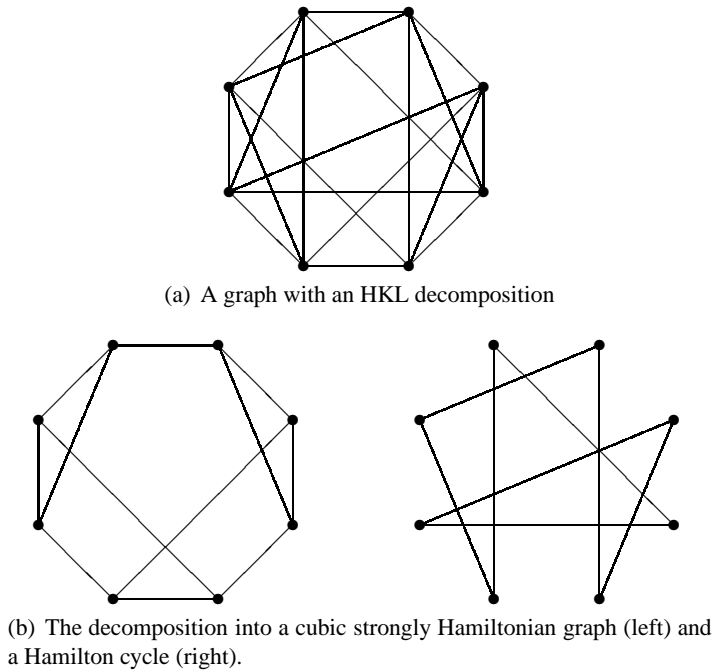

\end{example}

While it is unknown whether the complete graph $K_{2m}$ has a perfect 1-factorisation, it is possible to show that it has a 1-factorisation $\mathcal{F}$ such that the auxiliary digraph $\mathcal{H}(K_{2m},\mathcal{F})$ contains 2-factors, as demonstrated in the following lemma.

\begin{lemma} \label{lemma:CompleteHKL}
The complete graph $K_{2m}$ on an even number of vertices $2m$ has a 1-factorisation, $\mathcal{F}$ such that $\mathcal{H}(K_{2m},\mathcal{F})$ contains at least the arcs $(F_i,F_{i+1})$ and $(F_i,F_{i+2})$ and their reverses, for $0 \leq i \leq 2m-2$.
\end{lemma}
 
\begin{proof} 
There is a well-known 1-factorisation of $K_{2m}$, known as ${\rm GK}_{2m}$ (see \cite{Andersen07}), which is as follows.  Regard the vertices of $K_{2m}$ as $\mathbb{Z}_{2m-1}\cup \{\infty\}$, and define the {\em starter} 1-factor to be
\[ F_0 = \left\{ \{i,-i\} \,\, : \,\, i=1,\ldots,m-1 \right\} \cup \left\{ \{0,\infty\right\} \}. \]
The remaining 1-factors $F_i$ where $i=1,\ldots,2m-2$) are obtained by adding $i$ to each entry of $F_0$ (modulo $2m-1$, and with $\infty+i = \infty$).

Now, it is straightforward to show that the union of any two consecutive 1-factors, \mbox{$F_i \cup F_{i+1}$,} is a Hamilton cycle.  This establishes the existence of the arcs of first kind.  To show that  $(F_{i},F_{i+2})$ is an arc, we will show this is true for $i=-1$ and the result follows from the cyclic automorphism of the 1-factorisation. By considering which vertices are adjacent to $0$, we have the following path inside of $(F_{-1},F_{1}$):
\[ \begin{array}{ccccccccccccc} \pm 1= \mp (2m-2) & \cdots & 8 & -6 & 4 & -2 & 0 & 2 & -4 & 6 & -8 & \cdots & \pm (2m-2)= \mp 1, \end{array}\]
which can easily be seen to contain all of $\mathbb{Z}_{2m-1}$.  (Since we are working modulo $2m-1$, which is odd, the even negative numbers account for the positive odd numbers.)  Since $1$ and $-1$ are both adjacent to $\infty$, this does indeed yield a Hamilton cycle.
\end{proof}

The auxiliary digraph $\mathcal{H}(K_{2m},\mathcal{F})$ contains a Hamilton cycle and $K_{2m}$ has a sequentially perfect 1-factorisation (see Dinitz, Dukes and Stinson \cite{MR2134164}), and it also possesses an HKL decomposition; $K=F_{-1} \cup F_0 \cup F_1$ gives a cubic strongly Hamiltonian graph, while
\[ \mathcal{D} = \{ F_2\cup F_3, F_4 \cup F_5, \ldots, F_{2m-4}\cup F_{2m-3} \} \] forms a Hamiltonian decomposition of $K_{2m}\setminus K$.

The following construction gives UBBs for graph $G$ with a 1-factorisation for which the auxiliary digraph $\mathcal{H}(G,\mathcal{F})$ contains a directed 2-factor.
It is similar to that in the previous subsection for Hamilton-decomposable graphs.

\begin{construction} \label{constr:HKL}
Let G be $k$-regular and have a 1-factorisation $\mathcal{F} = \{F_0,F_1,\ldots,F_{k-1} \}$ for which $\mathcal{H}(G,\mathcal{F})$ contains a directed 2-factor. Suppose that in that 2-factor, the head of the arc whose tail is 1-factor $F$ is denoted $h(F)$. 

For each edge, $e \in F \subset G$, let $P_e$ be the path $F \cup h(F) \setminus \{e\}$. Let $t = k-1$.  We claim that this set of paths,
\[
\{P_e \mid e \in G\},
\] is a minimal $t$-UBB for $G$ with the number of bases equal to the number of edges of $G$.
\end{construction}

\begin{proofcorr}
Let $A$ be an arbitrary $t$-subset of edges of $G$, which we think of as ``bad'' edges we wish to avoid.  Since $|\mathcal{F}| = k > t$ there is at least one 1-factor that contains no bad edges.  For any 1-factor that contains no bad edges, if $h^{-1}(F)$ contains one or zero bad edges then we clearly have at least one base which is disjoint from $A$.  Thus if our set is not a $t$-UBB we must have at least two bad edges in $h^{-1}(F)$ whenever $F$ contains zero bad edges.  Let $z$ be the number of 1-factors that contain no bad edge.  There are therefore at least $z$ 1-factors that contain at least two bad edges and all the remaining 1-factors contain at least one bad edge.  This gives at least
\[
0\cdot z + 1\cdot (k-2z) + 2\cdot z = k
\]
bad edges which is a contradiction. \qed
\end{proofcorr}
Note that this proves that such a $k$-regular graph must have edge connectivity at least $k$.  Since a $k$-regular graph must have edge connectivity at most $k$ we obtain the following lemma.

\begin{lemma} \label{lemma:HKLEdgeConn}
Let $G$ be a graph of valency $k$ and with a 1-factorisation $\mathcal{F}$ for which $\mathcal{H}(G,\mathcal{F})$ contains a directed 2-factor.  Then the edge connectivity of $G$ is equal to $k$.
\end{lemma}

\begin{proofmin}
Let $T$ be a spanning tree from our UBB that we remove from the collection.  From Construction~\ref{constr:HKL}, we know there exists an edge $e\in F$ such that $T=(F \cup h(F))\setminus e$.  Then there exists a set of $t=k-1$ edges which is not uncoverable, consisting of $e$, and any edge from every 1-factor except $F$ and $h(F)$. \qed \end{proofmin}

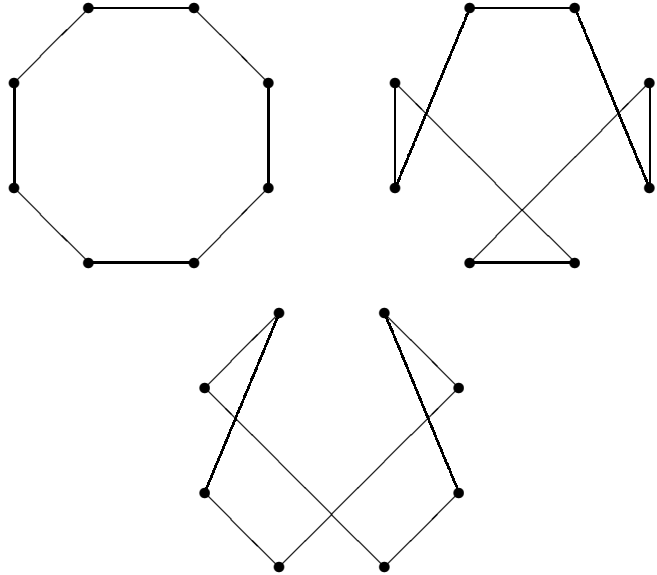
\begin{figure}[hbtp]
\setlength{\unitlength}{7mm}
\centering
\begin{picture}(6,5)

\put(2.000,0.000){\blob}
\put(0.586,1.414){\blob}
\put(0.586,3.414){\blob}
\put(2.000,4.828){\blob}
\put(4.000,4.828){\blob}
\put(5.414,3.414){\blob}
\put(5.414,1.414){\blob}
\put(4.000,0.000){\blob}

\put(0.586,3.414){\line(1,1){1.414}}
\put(2.000,4.828){\line(1,0){2.000}}
\put(4.000,4.828){\line(1,-1){1.414}}
\put(5.414,1.414){\line(0,1){2.000}}
\put(4.000,0.000){\line(1,1){1.414}}
\put(2.000,0.000){\line(1,0){2.000}}
\put(0.586,1.414){\line(1,-1){1.414}}
\put(0.586,1.414){\line(0,1){2.000}}


\end{picture}
\qquad
\begin{picture}(6,5)

\put(2.000,0.000){\blob}
\put(0.586,1.414){\blob}
\put(0.586,3.414){\blob}
\put(2.000,4.828){\blob}
\put(4.000,4.828){\blob}
\put(5.414,3.414){\blob}
\put(5.414,1.414){\blob}
\put(4.000,0.000){\blob}

\put(2.000,4.828){\line(1,0){2.000}}
\put(5.414,1.414){\line(0,1){2.000}}
\put(2.000,0.000){\line(1,0){2.000}}
\put(0.586,1.414){\line(0,1){2.000}}

\qbezier(0.586,1.414)(1.293,3.121)(2.000,4.828)
\qbezier(5.414,1.414)(4.707,3.121)(4.000,4.828)
\put(2.000,0.000){\line(1,1){3.414}}
\put(4.000,0.000){\line(-1,1){3.414}}

\end{picture}
\\[3ex]
\begin{picture}(6,5)

\put(2.000,0.000){\blob}
\put(0.586,1.414){\blob}
\put(0.586,3.414){\blob}
\put(2.000,4.828){\blob}
\put(4.000,4.828){\blob}
\put(5.414,3.414){\blob}
\put(5.414,1.414){\blob}
\put(4.000,0.000){\blob}

\put(0.586,3.414){\line(1,1){1.414}}
\put(4.000,4.828){\line(1,-1){1.414}}
\put(4.000,0.000){\line(1,1){1.414}}
\put(0.586,1.414){\line(1,-1){1.414}}

\qbezier(0.586,1.414)(1.293,3.121)(2.000,4.828)
\qbezier(5.414,1.414)(4.707,3.121)(4.000,4.828)
\put(2.000,0.000){\line(1,1){3.414}}
\put(4.000,0.000){\line(-1,1){3.414}}

\end{picture}
\caption{Three Hamilton cycles in a cubic strongly Hamiltonian graph. \label{figure:HKL2}}
\end{figure}

As in the previous two subsections, we don't claim that this construction is {\em minimum}.
\begin{example} \label{example:HKLUBB}
Consider the graph shown in Figure \ref{figure:HKL1}.  The cubic strongly Hamiltonian subgraph indicated in Figure \ref{figure:HKL1}(b) possesses three Hamilton cycles as shown in Figure \ref{figure:HKL2}.  In each cycle we take every other Hamilton path, as described in Construction~\ref{constr:HKL}, and include it in our UBB.  We also must take all spanning trees from the other Hamilton cycles in the HKL decomposition.  Our UBB contains $4+4+4+8=20 = |E(G)|$ spanning trees.  In contrast, the Sch\"{o}nheim lower bound for a \mbox{$(20,13,4)$} covering design is 11, while the best-known such covering design has size 16 \cite{lajolla}.
\end{example}

We note again that the number of spanning trees in the $t$-UBB is precisely the number of edges of $G$; so as in the previous subsections $|E(G)|$ is an upper bound on the minimum size of a $t$-UBB.

\subsection{Wheels} \label{subsect:wheel}
The {\em wheel} $W_n$ is the graph on $n+1$ vertices, formed from a cycle of length $n$ (the ``rim'') and an additional vertex (the ``hub'') adjacent to the $n$ others (by means of the ``spokes'').  Since $W_n$ has minimum degree 3, we have $\lambda(W_n) \leq 3$; furthermore, it is easy to see that the removal of any two edges leaves $W_n$ connected.  Hence $\lambda(W_n)=3$, and so we wish to construct a 2-UBB for these graphs.

Since $W_n$ has $2n$ edges, and a spanning tree for it has $n$ edges, the complements of the spanning trees in a 2-UBB for $W_n$ will form the blocks of a $(2n,n,2)$-covering design.  Now, the Sch\"onheim bound for covering designs (see \cite{MillsMullin92}) gives a lower bound of 6 (independent of $n$) for the size of a $(2n,n,2)$-covering design, and a construction due to Stanton, Kalbfleisch and Mullin (see Mills \cite[Theorem 3.2]{Mills79}) shows that this bound can always be attained.  Thus there is a lower bound of 6 for the size of a minimum 2-UBB for $W_n$.  Our construction was inspired by theirs.

We label the vertices of $W_n$ as follows: the vertices on the rim are $v_0,\ldots,v_{n-1}$ (with subscripts modulo $n$), while the hub is labelled as $v_\infty$.  Also, we label the edges as follows: those on the rim are labelled as $r_i$ (joining $v_i$ to $v_{i+1}$, while the spokes are labelled $s_i$ (joining $v_i$ to $v_\infty$).  We consider the cases where $n$ is even and odd separately.

First, suppose $n$ is even.  Partition the edges of $W_n$ into four sets, each of size $\frac{1}{2}n$, as follows:
\begin{eqnarray*}
A & = & \{ s_1, s_3, \ldots, s_{n-1} \} \\
B & = & \{ s_2, s_4, \ldots, s_{n-2}, r_0 \} \\
C & = & \{ r_1, r_3, \ldots, r_{n-1} \} \\
D & = & \{ r_2, r_4, \ldots, r_{n-2}, s_0 \}
\end{eqnarray*}
Let $\mathcal{U} =\{ A \cup B, A \cup C, A\cup D, B \cup C, B \cup D, C \cup D \}$.  We notice that each member of $\mathcal{U}$ forms a spanning tree for $W_n$.  Also, it is straightforward to see that any pair of edges of $W_n$ must be disjoint from one member of $\mathcal{U}$.  So $\mathcal{U}$ is a 2-UBB for $W_n$.  (We remark that since $\mathcal{U}$ is self-complementary, it is also a 2-covering-by-bases.)

Now, suppose $n$ is odd.  This time, we consider the following subsets of the edges of $W_n$:
\begin{eqnarray*}
A & = & \{ s_1, s_3, \ldots, s_{n-2} \} \\
B & = & \{ s_2, s_4, \ldots, s_{n-1} \} \\
C & = & \{ r_1, r_3, \ldots, r_{n-2} \} \\
D & = & \{ r_2, r_4, \ldots, r_{n-3} \}
\end{eqnarray*}
Then form the set $\mathcal{V}$ as follows:
\[ \mathcal{V} = \left\{
\begin{array}{l}
A\cup B\cup \{r_0\}, \\
A \cup C \cup \{r_{n-1}\}, \\
A \cup D \cup \{r_{n-1}\} \cup \{s_0\}, \\
B \cup C \cup \{s_0\}, \\
B \cup D \cup \{r_{n-1}\} \cup \{r_0\}, \\
C \cup D \cup \{r_0\} \cup \{s_0\} 
\end{array} \right\}. \]
Once again, we can verify that each member of $\mathcal{V}$ is a spanning tree for $W_n$.  Clearly any pair of edges from $A \cup B \cup C \cup D$ is disjoint from at least one member of $\mathcal{V}$ because there is some member that disjoint from the union of any two of $A$, $B$, $C$ or $D$.  A pair of edges, one from $A \cup B \cup C \cup D$ and the other from $\{r_0,r_{n-1}, s_0\}$ can be avoided because for any $S \in  \{A, B, C, D\}$ and any $e \in \{r_0,r_{n-1}, s_0\}$ there is an element of $\mathcal{V}$ which avoids $S \cup \{e\}$.  Finally it is easy to check that for any pair of edges from $\{r_0,r_{n-1}, s_0\}$, there exists some element of $\mathcal{V}$ disjoint from the pair. Thus  $\mathcal{V}$ must be a 2-UBB for $W_n$.  (This time, although it is not self-complementary, it is still a 2-covering-by-bases.) Since this construction meets the Sch\"{o}nheim bound, we are guaranteed that it definitely is a minimum UBB, unlike those in the previous subsections.  

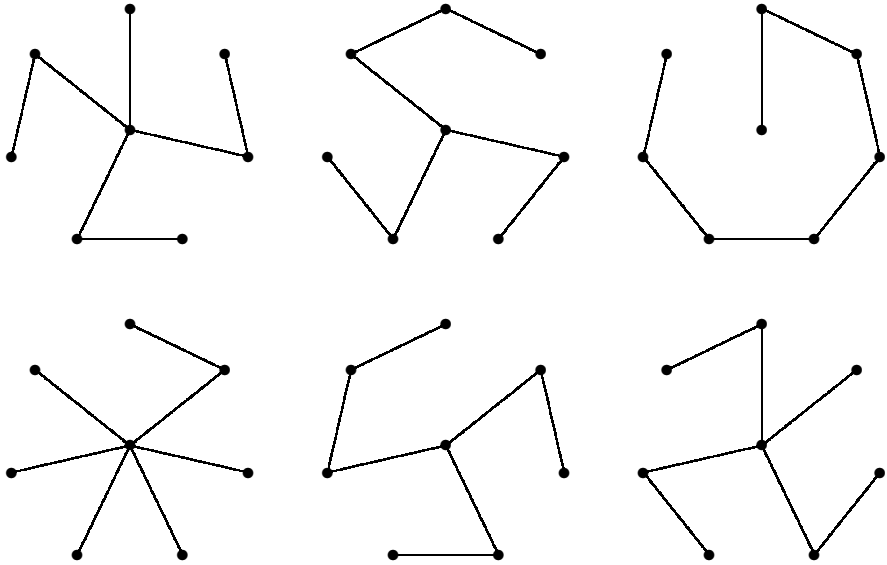
\begin{figure}[hbtp]
\centering
\setlength{\unitlength}{7mm}
\begin{picture}(18,11)


\put(2,0.5){\circle*{0.2}}
\put(4,0.5){\circle*{0.2}}
\put(5.2470,2.0637){\circle*{0.2}}
\put(4.8020,4.0136){\circle*{0.2}}
\put(3,4.8814){\circle*{0.2}}
\put(1.1980,4.0136){\circle*{0.2}}
\put(0.7530,2.0637){\circle*{0.2}}

\put(3,2.5765){\circle*{0.2}}

\qbezier(4.8020,4.0136)(3.9010,4.4475)(3,4.8814)      

\qbezier(3,2.5765)(2.0990,3.2951)(1.1980,4.0136) 
\qbezier(3,2.5765)(1.8765,2.3201)(0.7530,2.0637) 
\qbezier(3,2.5765)(2.5,1.5383)(2,0.5)            
\qbezier(3,2.5765)(3.5,1.5383)(4,0.5)            
\qbezier(3,2.5765)(4.1235,2.3201)(5.2470,2.0637) 
\qbezier(3,2.5765)(3.9010,3.2951)(4.8020,4.0136) 

\put(8,0.5){\circle*{0.2}}
\put(10,0.5){\circle*{0.2}}
\put(11.2470,2.0637){\circle*{0.2}}
\put(10.8020,4.0136){\circle*{0.2}}
\put(9,4.8814){\circle*{0.2}}
\put(7.1980,4.0136){\circle*{0.2}}
\put(6.7530,2.0637){\circle*{0.2}}

\put(9,2.5765){\circle*{0.2}}

\qbezier(9,4.8814)(8.0990,4.4475)(7.1980,4.0136)      
\qbezier(7.1980,4.0136)(6.9755,3.0387)(6.7530,2.0637) 
\put(8,0.5){\line(1,0){2}}                            
\qbezier(11.2470,2.0637)(11.0245,3.0387)(10.8020,4.0136) 

\qbezier(9,2.5765)(7.8765,2.3201)(6.7530,2.0637) 
\qbezier(9,2.5765)(9.5,1.5383)(10,0.5)            
\qbezier(9,2.5765)(9.9010,3.2951)(10.8020,4.0136) 

\put(14,0.5){\circle*{0.2}}
\put(16,0.5){\circle*{0.2}}
\put(17.2470,2.0637){\circle*{0.2}}
\put(16.8020,4.0136){\circle*{0.2}}
\put(15,4.8814){\circle*{0.2}}
\put(13.1980,4.0136){\circle*{0.2}}
\put(12.7530,2.0637){\circle*{0.2}}

\put(15,2.5765){\circle*{0.2}}

\qbezier(15,4.8814)(14.0990,4.4475)(13.1980,4.0136)      
\qbezier(12.7530,2.0637)(13.3765,1.2818)(14,0.5)         
\qbezier(16,0.5)(16.6235,1.2818)(17.2470,2.0637)         

\qbezier(15,2.5765)(13.8765,2.3201)(12.7530,2.0637) 
\qbezier(15,2.5765)(15.5,1.5383)(16,0.5)            
\qbezier(15,2.5765)(15.9010,3.2951)(16.8020,4.0136) 
\put(15,2.5765){\line(0,1){2.3049}}                 

\put(2,6.5){\circle*{0.2}}
\put(4,6.5){\circle*{0.2}}
\put(5.2470,8.0637){\circle*{0.2}}
\put(4.8020,10.0136){\circle*{0.2}}
\put(3,10.8814){\circle*{0.2}}
\put(1.1980,10.0136){\circle*{0.2}}
\put(0.7530,8.0637){\circle*{0.2}}

\put(3,8.5765){\circle*{0.2}}

\qbezier(1.1980,10.0136)(0.9755,9.0387)(0.7530,8.0637) 
\put(2,6.5){\line(1,0){2}}                            
\qbezier(5.2470,8.0637)(5.0245,9.0387)(4.8020,10.0136) 

\qbezier(3,8.5765)(2.0990,9.2951)(1.1980,10.0136) 
\qbezier(3,8.5765)(2.5,7.5383)(2,6.5)            
\qbezier(3,8.5765)(4.1235,8.3201)(5.2470,8.0637) 
\put(3,8.5765){\line(0,1){2.3049}}               

\put(8,6.5){\circle*{0.2}}
\put(10,6.5){\circle*{0.2}}
\put(11.2470,8.0637){\circle*{0.2}}
\put(10.8020,10.0136){\circle*{0.2}}
\put(9,10.8814){\circle*{0.2}}
\put(7.1980,10.0136){\circle*{0.2}}
\put(6.7530,8.0637){\circle*{0.2}}

\put(9,8.5765){\circle*{0.2}}

\qbezier(9,10.8814)(8.0990,10.4475)(7.1980,10.0136)      
\qbezier(6.7530,8.0637)(7.3765,7.2818)(8,6.5)         
\qbezier(10,6.5)(10.6235,7.2818)(11.2470,8.0637)         
\qbezier(10.8020,10.0136)(9.9010,10.4475)(9,10.8814)      

\qbezier(9,8.5765)(8.0990,9.2951)(7.1980,10.0136) 
\qbezier(9,8.5765)(8.5,7.5383)(8,6.5)            
\qbezier(9,8.5765)(10.1235,8.3201)(11.2470,8.0637) 

\put(14,6.5){\circle*{0.2}}
\put(16,6.5){\circle*{0.2}}
\put(17.2470,8.0637){\circle*{0.2}}
\put(16.8020,10.0136){\circle*{0.2}}
\put(15,10.8814){\circle*{0.2}}
\put(13.1980,10.0136){\circle*{0.2}}
\put(12.7530,8.0637){\circle*{0.2}}

\put(15,8.5765){\circle*{0.2}}

\qbezier(13.1980,10.0136)(12.9755,9.0387)(12.7530,8.0637) 
\qbezier(12.7530,8.0637)(13.3765,7.2818)(14,6.5)         
\put(14,6.5){\line(1,0){2}}                              
\qbezier(16,6.5)(16.6235,7.2818)(17.2470,8.0637)         
\qbezier(17.2470,8.0637)(17.0245,9.0387)(16.8020,10.0136) 
\qbezier(16.8020,10.0136)(15.9010,10.4475)(15,10.8814)      

\put(15,8.5765){\line(0,1){2.3049}}                 

\end{picture}
\caption{\label{figure:wheel} A 2-UBB for the wheel $W_7$.}
\end{figure}

\begin{example} \label{example:wheel}
Consider the wheel $W_7$.  The six spanning trees in a 2-UBB for $W_7$ are shown in Figure \ref{figure:wheel}.
\end{example}

\section{UBBs for random regular graphs}\label{random_regular}

Random regular graphs have been proposed as a model for networks: see Goerdt and Molloy \cite{Goerdt03}, for example.  In this model, each node in the network has a fixed number of neighbours, but the network is otherwise random.  Let $\mathcal{G}_{n,d}$ denote a random regular graph with $n$ vertices of degree $d$.  Now, if the degree $d$ of $\mathcal{G}_{n,d}$ is even, it is known that as $n \rightarrow \infty$, $\mathcal{G}_{n,d}$ possesses a Hamiltonian decomposition, asymptotically almost surely (i.e.~the probability of {\em not} having such a decomposition approaches zero); this is a result due to Kim and Wormald \cite{KimWormald01} (see also \cite{Wormald99}).  Furthermore, where the degree $d$ is odd, it is conjectured that as $n \rightarrow \infty$, $\mathcal{G}_{n,d}$ possesses a perfect 1-factorisation asymptotically almost surely.  (The reader should consult the survey by Wormald \cite{Wormald99} for full details.)  However, to show that $\mathcal{G}_{n,d}$ has an HKL decomposition, applying the principle of {\em contiguity arithmetic} (see \cite[Theorem 4.15(i)]{Wormald99}) it would suffice to show that $\mathcal{G}_{n,3}$ has a perfect 1-factorisation; this would likely be easier to prove.  Consequently, the constructions of UBBs in subsections \ref{subsect:hamdec} and \ref{subsect:HKL} are very relevant to networks, as they would allow (asymptotically, at least) our ``$t$-UBB protocol'' to be applied to this particular network model.

\section{A conjecture}

In all the cases we have considered, the number of spanning trees in the $t$-UBB was bounded above by the number of edges of $G$.  In the constructions in subsections \ref{subsect:Kmn}--\ref{subsect:HKL}, the spanning trees in our UBB were parameterised by the edge set, while our construction for wheels (in subsection \ref{subsect:wheel} had constant size.  Also, for max-STP graphs, our $t$-UBB consisted of $t+1$ edge-disjoint spanning trees, whereas the number of edges is clearly at least $(t+1)(n-1)$ (where $n$ is the number of vertices).

In view of this evidence, we make the following (possibly optimistic) conjecture.

\begin{conj}
Let $G$ be a graph with edge connectivity $\lambda(G)=k$.  Then there exists a $(k-1)$-UBB for $G$ with cardinality equal to $|E(G)|$.
\end{conj}

\section*{Acknowledgements}
The authors wish to thank NSERC and the Ontario Ministry of Research and Innovation for their financial support, and numerous colleagues for useful discussions.  The first author is a PIMS Postdoctoral Fellow.

\end{document}